\documentclass[twoside,12pt]{amsart}

\usepackage{amssymb}
\usepackage{amssymb}
\usepackage{latexsym}
\usepackage{verbatim}
\usepackage{graphics} 

\newtheorem{theo}{Theorem}
\newtheorem{prop}[theo]{Proposition}
\newtheorem{coro}[theo]{Corollary}
\newtheorem{lemme}[theo]{Lemma}

\newcommand{\RR}{{\mathbb R}}

\newcommand{\NN}{{\mathbb N}}
\newcommand{\ZZ}{\mathbb Z}
\newcommand{\QQ}{{\mathbb Q}}

\def\max{\hbox{\rm max}}

\def\mvide{{\epsilon }}

\def\gap{\hbox{\rm gap}}

\def\R{{\mathcal R}}
\def\S{{\mathcal S}}

\def\Sgood{{\mathcal S}_{\rm good}}
\def\Sprim{{\mathcal S}_{\rm prim}}

\def\Sconst{{\mathcal S}_{\rm const}}

\def\x{{\bf\sf x}}
\def\y{{\bf\sf y}}
\def\z{{\bf\sf z}}
\def\t{{\bf\sf t}}

\def\ind{{\bf 1}}

\def\sigmabar{\overline{\sigma}}
\def\Abar{\overline{A}}

\def\taubar{\overline{\tau}}

\def\tbar{\overline{\t}}

\def\aprime{a^{'}}

\begin{document}

\title{A theorem of Cobham for non-primitive substitutions}

\author{Fabien Durand}
\address{Universit\'{e} de Picardie Jules Verne,
Laboratoire Ami\'enois 
de Math\'emati\-ques Fondamentales  et Appliqu\'ees, 
CNRS-FRE 2270,
33 rue Saint Leu,
80039 Amiens Cedex 01, France.}

\email{fabien.durand@u-picardie.fr}


\maketitle

\section{Introduction.}

Given a subset $E$ of $\NN = \{ 0, 1,2,\cdots \}$ can we find an
elementary algorithm (i.e., a finite state automaton) which accepts the
elements of $E$ and rejects those that do not belong to $E$? In 1969
A. Cobham showed that the existence of such an algorithm deeply depends
on the numeration base. He stated \cite{Co1}: {\it Let $p$ and $q$ be
two multiplicatively independent integers (i.e., $p^k \not = q^l$ for
all integers $k,l>0$) greater than or equal to $2$. Let $E\subset
\NN$. The set $E$ is both $p$-recognizable and $q$-recognizable if and
only if $E$ is a finite union of arithmetic progressions.} What is
 now called the theorem of Cobham. We recall that a set $E\subset \NN$ is $p$-recognizable for some integer $p\geq 2$ if the language consisting of the expansions in base $p$ of the elements of $E$ is recognizable by a finite state automaton (see \cite{Ei}).

In 1972 Cobham gave an other partial answer to this question showing that not all sets are $p$-recognizable. He gave the following characterization: {\it The set $E\subset \NN $ is $p$-recognizable for some integer $p\geq 2$ if and only if the characteristic sequence $(x_n ; n\in \NN)$ of $E$ ($x_n = 1$ if $n\in E$ and $0$ otherwise) is generated by a substitution of length $p$,} where generated by a substitution of length $p$ means that it is the image by a letter to letter morphism of a fixed point of a substitution of length $p$.

We remark that $E$ is a finite union of arithmetic progressions if and
only if its characteristic sequence is ultimately periodic. Consequently
the theorem of Cobham can be formulated as follows (this is an
equivalent statement): {\it Let $p$ and $q$ be two multiplicatively
independent integers greater than or equal to $2$. Let $A$ be a finite alphabet and $\x\in A^{\NN}$. The sequence $\x$ is generated by both a substitution of length $p$ and a substitution of length $q$ if and only if $\x$ is ultimately periodic.}

To a substitution $\sigma$ is associated an integer square matrix $M\not
= 0$ which has non-negative entries. It is known (see \cite{LM} for
instance) that such a matrix has a real eigenvalue $\alpha$ which is
greater than or equal to the modulus of all others eigenvalues. It is usually called
the dominant eigenvalue of $M$. Let $\S$ be a set of substitutions. If $\x$ is the image by a letter to
letter morphism of a fixed point of $\sigma$ which belongs to $\S$ then we will say that $\x$
is $\alpha$-substitutive in $\S$. If $\S$ is the set of all
substitutions we will say that $\x$ is $\alpha$-substitutive. An easy computation shows that if $\sigma$ is
of length $p$ then $\alpha = p$. Furthermore if a sequence is generated
by a substitution of length $p$ then it is $p$-substitutive. Note that
the converse is not true. This suggests the following conjecture
formulated by G. Hansel.

\medskip

{\bf Conjecture.}
{\it
Let $\alpha$ and $\beta$ be two multiplicatively independent Perron numbers. Let $A$ be a finite alphabet. Let $\x$ be a sequence of $\in A^{\NN}$, the following are equivalent:
\begin{enumerate}
\item
$\x$ is both $\alpha$-substitutive and $\beta$-substitutive;
\item
$\x$ is ultimately periodic.
\end{enumerate}
}

\medskip

In this paper we prove that 2) implies 1) and, what is the main result of this paper, that this conjecture holds for a very large set of substitutions containing all known cases, we call it $\Sgood$. This set contains some non-primitive substitutions of non-constant length. More precisely for some sets $\S$ of substitutions, we prove

\begin{theo}
\label{maintheo}
Let $\alpha$ and $\beta$ be two multiplicatively independent Perron
 numbers. Let $A$ be a finite alphabet. A sequence $\x\in A^{\NN}$ is
 $\alpha$-substitutive in $\S$ and $\beta$-substitutive in $\S$ if and only if it is ultimately periodic.
\end{theo}

This result is true for $\Sconst$, the family of substitutions with constant length (this is the theorem of Cobham), and for $\Sprim$, the family of primitive substitutions \cite{Du2}. In \cite{Fa} and \cite{Du3} this result was proved for families of substitutions related to numeration systems. These families contain some non-primitive substitutions of non-constant length.

Much more results have been proved concerning generalizations of
Cobham's theorem to non-standard numeration systems \cite{BHMV1,BHMV2}.

Most of the proofs of Cobham's type results are divided into two parts. In the first part it is proven that the set $E\subset \NN$ is syndetic (the difference between two consecutive elements of $E$ is bounded) which corresponds to the fact that the letters of the characteristic sequence of $E$ appear with bounded gaps. In the second part the result is proven for such $E$. We will do the same. 

In Section \ref{preliminary} we recall some results concerning the length of the words $\sigma ^n (a)$ where $\sigma$ is a substitution on the alphabet $A$ and $a\in A$. These results have a key role in this paper.
In Section \ref{reciprocal} we prove that 2) implies 1). 
To prove the syndeticity of $E$ all proofs use the well-known fact that,
if $\alpha$ and $\beta$ are multiplicatively independent numbers
strictly greater than 1 then the set $\{ \alpha^n / \beta^m ; n,m \in
\ZZ \}$ is dense in $\RR^+$. Here we need more. We need the density in
$\RR^+$ of the set $\{ n^d\alpha^n / m^e\beta^m ; n,m \in \ZZ \}$, where
$d$ and $e$ are non-negative integers. We prove this result in Section
\ref{multindens} because we did not find it in the literature. We prove
in Section \ref{letterbg} that the letters with infinitely many
occurrences in $\x\in A^{\NN}$ appear with bounded gaps. This implies
the same result for words. In the last section we restrict ourself to
$\Sgood$, we recall some results obtained in \cite{Du3} and, using
return words, we conclude that $\x$ is ultimately periodic. More
precisely we prove that the conjecture is true for $\Sgood$.

\medskip

{\bf Words and sequences.} An {\it alphabet} $A$ is a finite set of
elements called {\it letters}. A {\it word} on $A$ is an element of the
free monoid generated by $A$, denoted by $A^*$. Let $x = x_0x_1 \cdots
x_{n-1}$ (with $x_i\in A$, $0\leq i\leq n-1$) be a word, its {\it
length} is $n$ and is denoted by $|x|$. The {\it empty word} is denoted by $\mvide$, $|\epsilon| = 0$. The set of non-empty words on  $A$ is denoted by $A^+$. The elements of $A^{\NN}$ are called {\it sequences}. If $\x=\x_0\x_1\cdots$ is a sequence (with $\x_i\in A$, $i\in \NN$), and $I=[k,l]$ an interval of $\NN$ we set $\x_I = \x_k \x_{k+1}\cdots \x_{l}$ and we say that $\x_{I}$ is a {\it factor} of $\x$. If $k = 0$, we say that $\x_{I}$ is a {\it prefix} of $\x$. The set of factors of length $n$ of $\x$ is written $L_n(\x)$ and the set of factors of $\x$, or {\it language} of $\x$, is noted $L(\x)$. The {\it occurrences} in $\x$ of a word $u$ are the integers $i$ such that $\x_{[i,i + |u| - 1]}= u$. When $\x$ is a word, we use the same terminology with similar definitions.

The sequence $\x$ is {\it ultimately periodic} if there
exist a word $u$ and a non-empty word $v$ such that $\x=uv^{\omega}$,
where $v^{\omega}= vvv\cdots $. Otherwise we say that
$\x$ is {\it non-periodic}. It is {\it periodic} if $u$ is the empty word. A sequence $\x$ is {\it uniformly recurrent} if for each factor $u$ the greatest difference of two successive occurrences of $u$ is bounded.

\medskip

{\bf Morphisms and matrices.} Let $A$ and $B$ be two alphabets. A {\it
morphism} $\tau$ is a map from $A$ to $B^*$. Such a map induces by
concatenation a morphism from $A^*$ to $B^*$. If $\tau (A)$ is included
in $B^+$, it induces a map from $A^{\NN}$ to $B^{\NN}$. These two maps
are also called $\tau$.

To a morphism $\tau$, from $A$ to $B^*$, is naturally associated the matrix $M_{\tau} = (m_{i,j})_{i\in  B  , j \in  A  }$ where $m_{i,j}$ is the number of occurrences of $i$ in the word $\tau(j)$.

Let $M$ be a square matrix, we call {\it dominant eigenvalue of $M$} an eigenvalue $r$ such that the modulus of all the other eigenvalues do not exceed the modulus of $r$. A  square matrix is called {\it primitive} if it has a power with positive coefficients. In this case the dominant eigenvalue is unique, positive and it is a simple root of the characteristic polynomial. This is Perron's Theorem.

A real number is a {\it Perron number} if it is an algebraic integer
that strictly dominates all its other albebraic conjugates. The
following result is well-known (see \cite{LM} for instance).

\begin{theo}
\label{lmp}
Let $\lambda $ be a real number. Then
\begin{enumerate}
\item
$\lambda $ is a Perron number if and only if it is
the dominant eigenvalue of a  primitive non-negative integral matrix.
\item
$\lambda$ is the spectral radius of a non-negative integral matrix if and
     only if $\lambda^p$ is a Perron number for some positive integer $p$.
\end{enumerate}
\end{theo}

\medskip

{\bf Substitutions and substitutive sequences.} In this paper a {\it substitution} is a morphism $\tau : A \rightarrow A^{*}$ such that for all letters of $A$ we have $\lim_{n\rightarrow +\infty} |\tau^n (a)| = +\infty$. Whenever the matrix associated to $\tau $ is primitive we say that $\tau$ is a {\it primitive substitution}.

A fixed point of $\tau $ is a sequence $\x=(\x_n ; n\in \NN )$ such that
$\tau (\x) = \x$. We say it is a {\it proper fixed point} if all
letters of $A$ have an occurrence in $\x$. We remark that all proper
fixed points of $\tau$ have the same language.

\medskip

{\bf Example.}
The substitution $\tau $ defined by $\tau (a)=aaab$, $\tau (b) = bc$ and $\tau (c) = b$
has two fixed points, one is starting with the letter $a$ and is proper
and the other one is starting with the letter $b$ and is not proper.

\medskip

If $\tau $ is a primitive
substitution then all its fixed points are proper and uniformly recurrent (for details
see \cite{Qu} for example).

Let $B$ be another alphabet, we say that a morphism $\phi$ from $A$ to
$B^*$ is a {\it letter to letter morphism} when $\phi (A)$ is a subset
of $B$. 
Let $\S$ be a set of substitutions and suppose that $\tau$ belongs to $\S$.
Then the sequence $\phi (\x)$ is called {\it substitutive in} $\S$.
We say $\phi (\x)$ is {\it substitutive} (resp. primitive substitutive)
if $\S$ is the set all substitutions (resp. the set of primitive substitutions).
 If $\x$ is a proper fixed point of $\tau$ and $\theta$ is the dominant eigenvalue of $\tau \in \S$ (i.e., the
 dominant eigenvalue of the matrix associated to $\tau$) then $\phi
 (\x)$ is called {\it $\theta$-substitutive in} $\S$; and we say {\it $\theta$-substitutive} (resp. primitive substitutive)
if $\S$ is the set all substitutions (resp. the set of primitive
substitutions). 

We point out that in the last example the fixed point
$\y$ of $\tau$ starting with the letter $b$ is also the fixed point of
the substitution $\sigma $ defined by $\sigma (b) = bc$ and $\sigma (c)
= b$. Moreover the dominant eigenvalue of $\tau $ is 3 and the dominant
eigenvalue of $\sigma $ is $(1+\sqrt{5})/2$. Hence in the definition of
``$\theta$-substitutive'' it is very important for $\x$ to be a proper
fixed point, otherwise the conjecture presented in the introduction
would not be true.

Clearly, if $\phi (\x)$ is $\theta$-substitutive then it
is $\theta^p$-substitutive for all $p\in \NN$. Consequently from Theorem
\ref{lmp} we can always suppose $\theta$ is a Perron number.

We define
$$
L(\tau )
=
\left\{
\tau^n (a)_{[i,j]} ; i,j\in \NN , i\leq j, n\in \NN , a\in A
\right\} .
$$
Let $\x$ be a fixed point of $\tau$. Then $L(\tau) = L(\x)$ if and only
if $\x$ is proper. If $\tau$ is primitive then for all its fixed points
$\x$ have the same language $L=L(\tau)$.

\section{Some preliminary lemmata.}

\label{preliminary}

This section and the first case of the proof of Proposition \ref{bgaps} is prompted by the ideas in \cite{Ha}.
 
In this section $\sigma$ will denote a substitution defined on the finite alphabet $A$, $\x$ one of its fixed points and $\Theta$ its dominant eigenvalue.

\begin{lemme}
\label{croissance}
There exists a unique partition $A_1, \cdots , A_l$ of $A$ such that for all $1\leq i \leq l$ and all $a\in A_i$
$$
\lim_{n\rightarrow +\infty} \frac{|\sigma^n (a) |}{c(a) n^{d(a)}\theta(a)^n} = 1
$$
where $\theta(a)$ is the dominant eigenvalue of $M$ restricted to $A_i$, $d(a)$ its Jordan order and $c(a)\in \RR$.
\end{lemme}
{\bf Proof.} See Theorem II.10.2 in \cite{SS}. \hfill $\Box$

\medskip

For all $a\in A$ we will call {\it growth type} of $a$ the couple $(d(a), \theta(a))$. If $(d,\alpha)$ and $(e,\beta)$ are two growth types we say that $(d,\alpha)$ is less than $(e,\beta)$ (or $(d,\alpha) < (e,\beta)$) whenever $\alpha < \beta$ or $\alpha = \beta$ and $d<e$. Consequently if the growth type of $a\in A$ is less then the growth type of $b\in A$ then $\lim_{n\rightarrow +\infty} |\sigma^n (a)| /|\sigma^n(b)| = 0$.

If the growth type of $a\in A$ is $(i , \theta )$ then there exists a letter $b$ with growth type $(i,\theta)$ having an occurrence in $\sigma (a)$.

We have $\Theta = \max \{ \theta(a) ; a\in A \}$. We set $D = \max \{ d(a) ; \theta(a) = \Theta, a\in A \}$ and $A_{max} = \{a\in A ; \theta (a) = \Theta , d(a) = D \}$. We will say that the letters of $A_{max}$ are {\it of maximal growth} and that $(D,\Theta)$ is the {\it growth type} of $\sigma$. 

For all letters $a\in A$, as $\lim_{n\rightarrow +\infty} |\sigma^n (a)| = +\infty$, it comes that $\theta (a) > 1$, or $\theta (a) = 1$ and $d(a) > 0$. Hence Lemma \ref{croissance} implies that there is no letter with growth type $(0,1)$. An important consequence of the following lemma is that in fact for all $a\in A$ we have $\theta (a) > 1$.

\begin{lemme}
\label{zerotheta}
If $(d,\theta)$ is the growth type of some letter then for all $i$
 belonging to $\{ 0, \cdots , d \}$ there exists a letter of growth type $(i,\theta)$ which appears infinitely often in $\x$.
\end{lemme}
{\bf Proof.} See Lemma III.7.10 in \cite{SS}.\hfill $\Box$

\medskip

We define
\begin{center}
\begin{tabular}{ccrcl}
$\lambda_{\sigma}$ & : & $A^*$               & $\rightarrow$ & $\RR$\\
          &   & $u_0\cdots u_{n-1}$ & $\mapsto $    & $\sum_{i=0}^{n-1} c(u_i)\ind_{A_{max}} (u_i)$.
\end{tabular}
\end{center}
From Lemma \ref{croissance} we deduce the following lemma.
\begin{lemme}
\label{lambda}
For all $u\in A^{*}$ we have $\lim_{n\rightarrow +\infty} |\sigma^n (u)|/n^D\Theta^n = \lambda_{\sigma} (u)$.
\end{lemme}

We say that $u\in A^{*}$ is of maximal growth if $\lambda_{\sigma} (u) \not = 0$.

\begin{lemme}
\label{etoile}
Let $a\in A$ which has infinitely many occurrences in $\x$. There exist a positive integer $p$, a word $u\in A^{*}$ of maximal growth and $v,w\in A^{*}$ such that for all $n\in \NN$ the word 
$$
\sigma^{pn} (u) \sigma^{p(n-1)} (v) \sigma^{p(n-2)} (v) \cdots \sigma^{p} (v) v w a 
$$
is a prefix of $\x$. Moreover we have
$$
\lim_{n\rightarrow + \infty}\frac{ | \sigma^{pn} (u) \sigma^{p(n-1)} (v) \sigma^{p(n-2)} (v) \cdots \sigma^{p} (v) v w a | }{\lambda_{\sigma} (u) (pn)^D \Theta^{pn} + \lambda_{\sigma} (v) \sum_{k=0}^{n-1}(pk)^D \Theta^{pk} } = 1 .
$$
\end{lemme}
{\bf Proof.} Let $a\in A$ be a letter that has infinitely many occurrences in $\x$. We set $a_0 = a$. There exists $a_1 \in A$ which has infinitely many occurrences in $\x$ and such that $a_0$ has an occurrence in $\sigma (a_1)$. In this way we can construct a sequence $(a_i ; i\in \NN)$ such that $a_0 = a$ and $a_i$ occurs in $\sigma ( a_{i+1})$, for all $i\in \NN$. 

There exist $i,j$ with $i<j$ such that $a_i = a_j = b$. It comes that $a$ occurs in $\sigma^i (b)$ and $b$ occurs in $\sigma^{j-i} (b)$. Hence there exist $u_1,u_2,v_1,v_2 \in A^{*}$ such that $\sigma^i (b) = u_1 a u_2$ and $\sigma^{j-i} (b) = v_1 b v_2$. We set $p = j-i$, $v = \sigma ^i (v_1)$ and $w = u_1$. There exists $u^{'}$ such that $u^{'}b$ is a prefix of $\x$. We remark that for all $n\in \NN$ the word $\sigma^n (u^{'}b)$ is a prefix of $\x$ too. We set $u = \sigma^i (u^{'})$. We have $\sigma^{p} (u^{'}b) = \sigma^{p} (u^{'}) v_1 b v_2$. Consequently for all $n\in \NN$
$$
\sigma^{pn} (u^{'}) \sigma^{p(n-1)} (v_1) \sigma^{p(n-2)} (v_1) \cdots \sigma^{p} (v_1) v_1 b 
$$ 
is a prefix of $\sigma^{np} (u^{'}b)$. Then 
$$
\sigma^{pn} (u) \sigma^{p(n-1)} (v) \sigma^{p(n-2)} (v) \cdots \sigma^{p} (v) v w a 
$$
is a prefix of $\sigma^{np+i} (u^{'}b)$ and consequently of $\x$, for all $n\in \NN$. The last part of the lemma follows from Lemma \ref{lambda}. \hfill $\Box$

\section{Assertion 2) implies Assertion 1) in the conjecture.}

\label{reciprocal}

In this section we prove the following proposition. It it is the
``easy'' part of the conjecture, namely Assertion 2) implies Assertion
1). The first part of the proof is an adaptation of the proof of
Proposition 3.1 in \cite{Du1} and the second part is inspired by the
substitutions introduced in Section V.4 and Section V.5 of \cite{Qu}.

\begin{prop}
\label{recip}
Let $\x$ be a sequence on a finite alphabet and $\alpha$ a Perron number. If $\x$ is periodic (resp. ultimately periodic) then it is $\alpha$-substitutive primitive (resp. $\alpha$-substitutive).
\end{prop}
{\bf Proof.} Let $\x$ be a periodic sequence with period $p$. Hence we can suppose that $A = \{ 1,\cdots , p \}$ and $\x = (1\cdots p )^{\omega}$. Let $M$ be a primitive matrix whose dominant eigenvalue is $\alpha$ and $\sigma : B \rightarrow B^*$ a primitive substitution whose matrix is $M$. Let $\y$ be one of its fixed points. In the sequel we construct, using $\sigma$, a new substitution $\tau$ with dominant eigenvalue $\alpha$, together with a fixed point $\z=\tau (\z)$, and a letter to letter morphism $\phi$ such that $\phi ( \z ) = \x$. We define the alphabet

$$
D = \left\{ (b,i) \ ; b\in B \ , \ 1\leq i \leq p \right\} ,
$$

the morphism $\psi : B \rightarrow D^{*}$ and the substitution $\tau : D\rightarrow D^{*}$ by 
\[
\begin{array}{lll}
\psi(b) = (b,1)\cdots (b,p) & {\rm and} & \tau ((b,i)) = (\psi (\sigma (b)))_{[(i-1)|\sigma(b)| , i|\sigma (b)|-1]} ,
\end{array}
\]
for all $(b,i)\in D$. The substitution $\tau $ is well defined because
$|\psi (  \sigma (b))| = p|\sigma(b)|$. Moreover, these morphisms are
such that $\tau \circ \psi = \psi \circ \sigma$. Hence the substitution
$\tau$ is primitive. The sequence $\z = \psi (\y)$ is a fixed point of $\tau$ and (using Perron theorem and the fact that $M_{\tau} M_{\psi} = M_{\psi} M_{\sigma}$) its dominant eigenvalue is $\alpha$.

Let $\phi : D \rightarrow A$ be the letter to letter morphism defined by $\phi ((b,i)) = i$. It is easy to see that $\phi (\z) = \x$. It follows that $\x$ is $\alpha$-substitutive.

Suppose now that $\x$ is ultimately periodic. Then there exist two non-empty words $u$ and $v$ such that $\x = uv^{\omega}$. From what precedes we know that there exist a substitution $\tau : D\rightarrow D^{*}$, a fixed point $\z = \tau (\z)$ and a letter to letter morphism $\phi : D \rightarrow A$ such that $\phi ( \z ) = v^{\omega}$. Let $E^{'} = \{ a_1 , a_2 , \cdots , a_{|u|} \}$ be an alphabet, with $|u|$ letters, disjoint from $D$ and consider the sequence $\t = a_1a_2\cdots a_{|u|} \z \in (E^{'} \cup D)^{\NN} = F^{\NN}$. It suffices to prove that $\t$ is $\alpha$-substitutive. We extend $\tau$ to $F$ setting $\tau (a_i) = a_i$, $1\leq i \leq |u|$. Let $G$ be the alphabet of the words of length $|u|+1$ of $\t$, that is to say
$$
G = \left\{ (\t_n \t_{n+1} \cdots \t_{n+|u|}) ; n\in \NN \right\} \ \ {\rm where} \ \ \t = \t_0 \t_1 \cdots .
$$
The sequence $\tbar = (\t_0 \t_{1} \cdots \t_{|u|})(\t_1 \t_{2} \cdots
\t_{|u|+1}) \cdots (\t_n \t_{n+1} \cdots \t_{n+|u|}) \cdots \in G^{\NN}$ is a fixed
point of the substitution $\zeta : G \rightarrow G^{*}$ we define as
follows. Let $(l_0l_1 \cdots l_{|u|-1} a)$ be an element of $G$. Let
$s_0 s_1 \cdots s_{|u|-1}$ be the suffix of length $|u|$ of the word $\tau ( l_0l_1 \cdots l_{|u|-1})$.

If $|\tau (a)| \leq |u|$, we set 
$$ 
\zeta ((l_0l_1 \cdots l_{|u|-1} a )) 
$$
$$ 
=
(s_{[0,|u|-1]} \tau (a)_0)
(s_{[1 , |u|-1]} \tau (a)_{[0,1]})
\cdots
(s_{[|\tau (a)|-1 , |u|-1]} \tau (a)_{[0,|\tau (a)|-1]}) ,
$$
otherwise 
$$
\zeta ((l_0l_1 \cdots l_{|u|-1} a )) 
$$
$$
=(s_{[0,|u|-1]} \tau (a)_0)
\cdots
(s_{|u|-1} \tau (a)_{[0,|u|-1]})
(\tau (a)_{[0,|u|]})
\cdots
(\tau (a)_{[|\tau (a)| - |u| - 1,|\tau (a)|-1]}) ,
$$

By induction we can prove that for all $n\in \NN$ we have 
$$
\zeta^n (( \t_0 \t_1 \cdots \t_{|u|})) 
$$
$$
= ( \t_0 \t_1 \cdots \t_{|u|})( \t_1 \t_2 \cdots \t_{|u|+1}) \cdots ( \t_{|\tau^n (\t_{|u|})| -1}  \cdots 
\t_{|\tau^n (\t_{|u|})| + |u| - 1 } ) .
$$
Consequently $\tbar$ is a fixed point of $\zeta$ and $\rho (\tbar ) = \t$ where $\rho : G \rightarrow F$ is defined by 
$$
\rho (( r_0 r_1 \cdots r_{|u|}  )) = r_{0}  .
$$ 
Moreover we remark that for all $n\in \NN$ we have
$$
| \zeta^n (( r_0 r_1 \cdots r_{|u|}  )) | = | \tau^{n} ( r_{|u|} ) | .
$$
From this and Lemma \ref{croissance} it comes that for all $( r_0 r_1 \cdots r_{|u|}  ) \in D$ we have 
$$
\lim_{n\rightarrow +\infty} \frac{| \zeta^{n+1} (( r_0 r_1 \cdots r_{|u|}  )) |}{| \zeta^n (( r_0 r_1 \cdots r_{|u|}  )) |} = \alpha .
$$
Hence $\alpha $ is the dominant eigenvalue of $\zeta$ and $\t$ is $\alpha$-substitutive. 
\hfill $\Box$

\medskip

{\bf Example.} Let $\x = (12)^{\omega}$ and $\alpha = (1+\sqrt{5})/2$. It is the dominant eigenvalue of the substitution $\sigma : A= \{ a,b \} \rightarrow A^*$ given by $\sigma (a) = ab$ and $\sigma (b) = a$. We have $D = \{ (a,1), (a,2),(b,1),(b,2)  \}$ and the substitution $\tau : D\rightarrow D^*$ defined in the previous proof is given by 
$$
\begin{array}{lll}
\tau ((a,1)) = (a,1)(a,2) , &           & \tau ((a,2)) = (b,1)(b,2) ,\\ 
\tau ((b,1)) = (a,1)        & {\rm and} & \tau ((b,2)) = (a,2) .
\end{array}
$$

\medskip

{\bf Example.} Let $c$ be a letter and $\x = c (12)^{\omega}$. We take the notations of the previous example and for convenience we set $A=(a,1)$, $B=(a,2)$, $C=(b,1)$ and $D=(b,2)$. The substitution $\zeta : G \rightarrow G^*$, where $G = \{ (cA), (AB),(BC),(CD),(DA),(BA) \}$, defined in the previous proof is given by
\begin{center}
\begin{tabular}{ll}
$\zeta ((cA)) = ((cA))((AB))$, & $\zeta ((AB)) = ((BC))((CD))$,\\
$\zeta ((BC)) = ((DA))$,       & $\zeta ((CD)) = ((AB))$,\\ 
$\zeta ((DA)) = ((BA))((AB))$, & $\zeta ((BA)) = ((DA))((AB))$.
\end{tabular}
\end{center}
Let $\tbar$ be the fixed point of $\zeta $ whose first letter is $(cA)$. Let $\phi : G \rightarrow \{ c,1,2 \}$ be the letter to letter morphism given by
$$
\begin{array}{lll}
\phi ((cA)) = c, & \phi ((AB)) = 1, & \phi ((BC)) = 2,\\ 
\phi ((CD)) = 1, & \phi ((DA)) = 2, & \phi ((BA)) = 2.
\end{array}
$$
We have $\phi (\tbar) = c(12)^{\omega} = \x$.

\medskip

Using Proposition \ref{recip} we obtain a slight improvement of the main results of respectively \cite{Du2} and \cite{Du3}. More precisely:

\begin{theo}
Let $\alpha$ and $\beta$ be two multiplicatively independent Perron numbers. Let $\x$ be a sequence on a finite alphabet. The sequence $\x$ is both $\alpha$-substitutive primitive and $\beta$-substitutive primitive if and only if it is periodic.
\end{theo}
\begin{theo}
Let $U$ and $V$ be two Bertrand numeration systems, $\alpha$ and $\beta$ be two multiplicatively independent $\beta$-numbers such that $L(U) = L(\alpha)$ and $L(V) = L(\beta)$. Let $E$ be a subset of $\NN$. The set $E$ is both $U$-recognizable and $V$-recognizable if and only if it is a finite union of arithmetic progressions. {\rm (see \cite{Du3} for the terminology)}
\end{theo}

\section{Multiplicative independence and density.}

\label{multindens}

This section is devoted to the proof of the following proposition.
\begin{prop}
\label{denselog}
Let $\alpha$ and $\gamma$ be two rationally independent positive numbers (i.e., $\alpha / \beta \not \in \QQ$). Let $d$ and $e$ be two non-negative integers. Then the set
$$
\left\{
n\alpha +d\log n -m\beta - e\log m  ; \ n,m \in \NN
\right\}
$$
is dense in $\RR$.
\end{prop}
The following straightforward corollary will be essential in the next section.
\begin{coro}
\label{densite}
Let $\alpha$ and $\beta$ be two multiplicatively independent positi\-ve real numbers. Let $d$ and $e$ be two non-negative integers. Then the set
$$
\left\{
\frac{n^d\alpha^n}{m^e\beta^m} ; n,m \in \NN
\right\}
$$
is dense in $\RR^+$. 
\end{coro}
These two results are well-known for $d=e=0$ (see \cite{HW} for example). We need the following lemma to prove Proposition \ref{denselog}.
\begin{lemme}
\label{lemmetech}
Let $\beta <\alpha $ be two rationally independent numbers. Then for all
 $\epsilon > 0$ and all $N\in \NN$ there exist $m,n$, with $m\geq n \geq N$, such that $0<n\alpha - m\beta < \epsilon$.
\end{lemme}
{\bf Proof.} The proof is left to the reader.\hfill $\Box$

\medskip

{\bf Proof of Proposition \ref{denselog}.} Let $l\in \RR$ and $\epsilon > 0$, we have to find $N,M\in \NN$ such that $|N\alpha +d\log N -M\beta - e\log M - l| < \epsilon$. The proof is divided into several cases.

\medskip

First case: $\alpha > \beta$, $e=d$ and $l\geq d\log(\frac{\beta}{\alpha})$.

\medskip

From Lemma \ref{lemmetech} there exist two integers $0<n<m$ such that $0<n\alpha - m\beta < \frac{\epsilon}{2}$ and $d\log (1+\frac{\epsilon}{m\beta})\leq \frac{\epsilon}{2}$. Hence we have
\begin{equation}
\label{eq0}
d\log (\frac{\beta}{\alpha}) 
< 
d\log (n) - e\log (m) 
<
 d \log (\frac{\beta}{\alpha}) + d \log (1+\frac{\epsilon}{m\beta}).
\end{equation}
Then $n\alpha - m\beta + d(\log n - \log m ) < l+\epsilon$. We consider $f: \NN \rightarrow \RR$ defined by
$$
f(k) = k(n\alpha - m\beta) - d(\log (km) - \log (kn) ).
$$
We have $f(1) < l + \epsilon$, $\lim_{k\rightarrow +\infty} f(k) = +\infty$ and $0< f(k+1) - f(k) = n\alpha - m\beta < \epsilon$. Hence there exists $k_0\in \NN$ such that $|f(k_0) - l| < \epsilon$, that is to say
$$
|N\alpha +d\log N -M\beta - e\log M - l| < \epsilon
$$
where $N = nk_0$ and $M = mk_0$.

\medskip

Second case:  $\alpha > \beta$, $e=d$ and $l< d\log(\frac{\beta}{\alpha})$.

\medskip

It suffices to take $n,m$ with $0<n<m$ such that $-\frac{\epsilon}{2} <n\alpha - m\beta < 0$ and $d\log (1+\frac{\epsilon}{m\beta})\leq \frac{\epsilon}{2}$, and the same method will give the result.

\medskip

Third case: $\alpha > \beta$ and $e>d$.

\medskip

Let $k_0\in \NN$ be such that $-\epsilon < (d-e)\log (1 + \frac{1}{k_0})
< 0$. If two integers $n,m$ with $0<n<m$ are such that $0<n\alpha - m\beta <\epsilon$ then we have
$$
(d-e)\log (m) + d\log (\frac{\beta}{\alpha}) 
< 
d\log (n) - e\log (m) 
$$
$$
<
(d-e)\log (m) + d \log (\frac{\beta}{\alpha}) + d \log (1+\frac{\epsilon}{m\beta}),
$$
which is negative for $m$ large enough. Hence from Lemma \ref{lemmetech}
it comes that there exist two integers $n,m$ with $0<n<m$ such that $0< n\alpha - m\beta <\epsilon$ and 
\begin{equation}
d\log (n) - e\log (m) \leq l - (k_0)\epsilon - (d-e)\log (k_0).
\end{equation}
We consider $f: \NN \rightarrow \RR$ defined by
$$
f(k) = k(n\alpha - m\beta) + d\log (kn) - e\log (km).
$$
We have
$$
f(k_0) \leq k_0 \epsilon + (d-e)\log (k_0) + d\log (n) - e\log (m) \leq l.
$$
Moreover $\lim_{k\rightarrow +\infty} f(k) = +\infty$ and for all $k\geq k_0$
$$
-\epsilon < f(k+1) - f(k) = n\alpha - m\beta + (d-e)\log(1+\frac{1}{k})< \epsilon .
$$
Hence there exists an integer $k_1 \geq k_0$ such that $|f(k_1) - l| < \epsilon$, that is to say
$$
|N\alpha +d\log N -M\beta - e\log M - l| < \epsilon
$$
where $N = nk_1$ and $M = mk_1$.

\medskip

Remaining cases: The same ideas achieve the proof.\hfill $\Box$

\section{The letters appear with bounded gaps.}
\label{letterbg}

Let $\alpha$ and $\beta$ be two multiplicatively independent Perron numbers. Let $\sigma$ and $\tau$ be two substitutions on the alphabets $A$ and $B$, with fixed points $\y$ and $\z$ and with growth types  $(d , \alpha)$ and $(e , \beta )$ respectively. Let $\phi : A\rightarrow C$ and $\psi : B\rightarrow C $ be two letter to letter morphisms such that $\phi (\y) = \psi (\z) = \x$.

This section is devoted to the proof of the following proposition.

\begin{prop}
\label{bgaps}
The letters of $C$ which have infinitely many occurrences appear in $\x$ with bounded gaps in $\x$.
\end{prop}
{\bf Proof:} We prove this proposition considering two cases. 

\medskip

Let $c\in C$ which has infinitely many occurrences. Let $X = \{ n \in
 \NN ; \x_n = c \}$ and $A^{'} = \{ a\in A ; \phi (a) = c \}$. Assume
 that the letter $c$ does not appear with bounded gaps. Then there exist
 $a\in A$ with infinitely many occurrences in $\y$ and a strictly
 increasing sequence $(p_n ; n\in \NN)$ of positive integers such that
 the letter $c$ does not appear in $\phi (\sigma^{p_n} (a) )$. Let $A^{''}$ be the set of such letters. We consider two cases.

\medskip

{\bf First case: There exists $a\in A^{''}$ of maximal growth.}

Let $u\in A^*$ such that $ua$ is a prefix of $\y$. Of course we can suppose that $u$ is non-empty. 

For all $n\in \NN$ we call $\Omega_n\subset A$ the set of letters appearing in $\sigma^{p_n} (a)$. There exist two distinct integers $n_1 < n_2$ such that $\Omega_{n_1} = \Omega_{n_2}$. Let $\Omega$ be the set of letters appearing in $\sigma^{p_{n_2} - p_{n_1}} ( \Omega_{n_1} )$. It is easy to show that $\Omega = \Omega_{n_1} = \Omega_{n_2}$.

Consequently the set of letters appearing in $\sigma^{p_{n_2} - p_{n_1}} (\Omega)$ is equal to $\Omega$ and for all $k\in \NN$ the set of letters appearing in $\sigma^{p_{n_1} + k ( p_{n_2} - p_{n_1} ) } (A)$ is equal to $\Omega$. We set $p = p_{n_1}$ and $g = p_{n_2} - p_{n_1}$. We remark that the letter $c$ does not appear in the word $\phi (\sigma^{p+kg} (a))$ and that $[|\sigma^{p+kg} (u)|, |\sigma^{p+kg} (ua)|[ \cap X = \emptyset$, for all $k\in \NN$.

There exists a letter $\aprime$ of maximal growth having an occurrence in $\sigma^p (a)$. We set $\sigma^p (a) = w\aprime w^{'}$. For all $k\in \NN$ we have $|\sigma^{p+kg} (ua)| \geq |\sigma^{kg} (\sigma^p (u)wa^{'})|$ and
\begin{equation}
\label{gap}
[|\sigma^{kg} (v)| , |\sigma^{kg} (vw\aprime)|[ \cap X = \emptyset
\end{equation}
where $v = \sigma^p (u)$. Because $\aprime$ is of maximal growth we have $\lambda_{\sigma} (v) < \lambda_{\sigma} (vw\aprime)$. Consequently there exists an $\epsilon > 0$ such that
$$
\lambda_{\sigma} (v) (1+\epsilon) < \lambda_{\sigma} ( vw\aprime) (1-\epsilon).
$$
From Lemma \ref{lambda} we obtain that there exists $k_0$ such that for all $k\geq k_0$ we have
\begin{equation}
\label{encadre}
\frac{|\sigma^{kg} (v)|}{(kg)^d \alpha^{kg}} 
< \lambda_{\sigma} (v) (1+\epsilon) 
< \lambda_{\sigma} (vw\aprime) (1-\epsilon)
< \frac{|\sigma^{kg} (vw\aprime)|}{(kg)^d \alpha^{kg}} 
.
\end{equation}
From Lemma \ref{etoile} applied to $\tau$ we have that there exist $s\in B^{*}$ of maximal growth, $t,t^{'} \in B^{*}$ and $h\in \NN^{*}$ such that for all $n\in \NN$ 
$$
\psi \left( \y_{[\tau^{hn} (s) \tau^{h(n-1)} (t) \cdots  \tau^{h} (t) t t^{'}]}\right) = c.
$$
From the second part of Lemma \ref{etoile} it comes that there exists $\gamma \in \RR$ such that 
$$
\lim_{n\rightarrow +\infty} \frac{| \tau^{hn} (s) \tau^{h(n-1)} (t) \cdots  \tau^{h} (t) t t^{'}|}{(nh)^e\beta^{hn}} = \gamma .
$$

From Corollary \ref{densite} it comes that there exist two strictly increasing sequences of integers, $(m_i ; i\in \NN)$ and $(n_i ; i\in \NN )$, and $l\in \RR$ such that 
$$
\frac{\gamma (m_ih)^e \beta^{m_ih}}{(n_ig)^d \alpha^{n_ig}} \longrightarrow_{i\rightarrow +\infty} 
l 
\in \
]\lambda_{\sigma} (v) (1+\epsilon) , \lambda_{\sigma} (vw\aprime) (1-\epsilon)[.
$$
Hence from Lemma \ref{lambda} we also have
$$
\frac{|\tau^{hm_i} (s) \tau^{h(m_i-1)} (t) \cdots  \tau^{h} (t) t t^{'}
|}{(n_ig)^d \alpha^{n_i g}}
$$
\begin{equation}
\label{limite}
=
\frac{| \tau^{hm_i} (s) \tau^{h(m_i-1)} (t) \cdots  \tau^{h} (t) t t^{'} |}{\gamma (m_ih)^e \beta^{m_ih}}
\frac{\gamma (m_ih)^e \beta^{m_ih}}{(n_ig)^d \alpha^{n_i g}}
\longrightarrow_{i\rightarrow +\infty} l.
\end{equation}
From (\ref{encadre}) and (\ref{limite}) there exists $i\in \NN$ such that
$$
|\sigma^{n_ig} (v)|
<
|\tau^{hm_i} (s) \tau^{h(m_i-1)} (t) \cdots  \tau^{h} (t) t t^{'} |
<
|\sigma^{n_ig} (vw\aprime)|,
$$
which means that $ |\tau^{hm_i} (s) \tau^{h(m_i-1)} (t) \cdots  \tau^{h}
(t) t t^{'} | $ belongs to $X$. This gives a contradiction with (\ref{gap}).

\medskip

{\bf Second case: No letter in $A^{''}$ has maximal growth.}

We define $B^{''}$ as $A^{''}$ but with respect to $\tau$ and $B$. We can suppose that no letter of $B^{''}$ has maximal growth. 

There exists a letter $a\in A^{''}$ (resp. $b\in B^{''}$) which has infinitely many occurrences in $\y$ (resp. $\z$) and with growth type $(d^{'},\alpha^{'})<(d,\alpha)$ (resp. $(e^{'},\beta^{'})<(e,\beta)$). We recall that $\alpha^{'}$ and $\beta^{'}$ are greater than $1$. Furthermore we can suppose that $(d^{'},\alpha^{'})$ (resp. $(e^{'},\beta^{'})$) is maximal with respect to $A^{''}$ (resp. $B^{''}$).

Let $w=w_0 \cdots w_n$ be a word belonging to $L(\y)$ (resp. $L(\z)$), we call $\gap (w)$ the largest integer $k$ such that there exists $0\leq i \leq n-k+1$ for which the letter $c$ does not appear in $\phi (w_i \cdots w_{i+k-1})$ (resp. in $\psi (w_i \cdots w_{i+k-1})$).

There exist infinitely many prefixes of $\y$ (resp. $\z$) of the type $u_1au_2a'$
(resp. $v_1bv_2b'$) fulfilling the conditions $\imath)$ and $\imath\imath)$ below:

\medskip

$\imath$) 
The growth type of $u_1\in A^*$ and $a'\in A$ (resp. $v_1\in B^*$ and $b'\in B$) is maximal.

$\imath\imath$) 
The words $u_2$ and $v_2$ do not contain a letter of maximal growth.

\medskip

It is easy to prove that there exists a constant $K^{'}$ such that $\gap(\tau^{n} (b')) \leq K^{'} n^{e'} {\beta^{'}}^{n}$ and $\gap(\sigma^{n} (a')) \leq K^{'} n^{d'} {\alpha^{'}}^{n}$ for all $n\in \NN$.
Due to Lemma \ref{croissance}, $\lim_{n\rightarrow +\infty }|\sigma^{n}
(a)|/n^{d'}{\alpha^{'}}^n$ and  $\lim_{n\rightarrow +\infty }|\tau^{n}
(b)|/n^{e'}{\beta^{'}}^n$ exist and are finite, we call them $\mu (a)$
and $\mu (b)$ respectively. 

Let $u_1au_2a'$ be a prefix of $\y$ fulfilling the conditions $\imath$)
and $\imath\imath $), then choose $v_1bv_2b'$ fulfilling the same
conditions and so that 
\begin{equation}
\label{correction}
\frac{K'}{\mu (a)}
\left(
\frac{2\lambda_{\sigma} (u_1)}{2\lambda_{\tau} (v_1) + \lambda_{\tau}(b')}
\right)^{\frac{\log \alpha'}{  \log \alpha}}
\left(
\frac{\log \beta}{\log \alpha}
\right)^{e \frac{\log (\beta')}{\log (\beta)} -e' }
\leq \frac{1}{3} .
\end{equation}

From Corollary \ref{densite}  there exist four strictly increasing sequences of integers $(m_i;i\in \NN)$, $(n_i;i\in\NN)$, $(p_i;i\in\NN)$ and $(q_i;i\in\NN)$ such that 

\medskip

\begin{equation}
\label{multind}
\begin{array}{lll}
\lim_{i\rightarrow +\infty}  \frac{n_i^d \alpha^{n_i}}{m_i^e\beta^{m_i}}  
 &=& 
\frac{2\lambda_{\tau} (v_1)}{2 \lambda_{\sigma}(u_1) + \lambda_{\sigma}(a')}
\ \
{\rm and} \\
 & & \\
\lim_{i\rightarrow +\infty}  \frac{p_i^e\beta^{p_i}}{q_i^d \alpha^{q_i}} 
 &=&
\frac{2\lambda_{\sigma} (u_1)}{2\lambda_{\tau} (v_1) + \lambda_{\tau}(b')}.
\end{array}
\end{equation}

\medskip

As a consequence of (\ref{multind}) we have 

\medskip

\begin{equation}
\label{nu} 
\lim_{i\rightarrow +\infty} n_i/m_i = \log (\beta)/ \log (\alpha) \ \  {\rm and} \ \ \lim_{i\rightarrow +\infty} p_i/q_i = \log (\alpha) / \log (\beta),
\end{equation}

\medskip

and there exists $i_0$ such that for all $i\geq i_0$ we have
$$
\frac{|\sigma^{n_i}(u_1au_2)|}{|\tau^{m_i}(v_1)|} 
\leq 1 \leq 
\frac{|\sigma^{n_i}(u_1au_2a')|}{|\tau^{m_i}(v_1b)|}
\ \ {\rm and}
$$
$$ 
\frac{|\tau^{p_i}(v_1bv_2)|}{|\sigma^{q_i}(u_1)|} 
\leq 1 \leq 
\frac{|\tau^{p_i}(v_1bv_2b')|}{|\sigma^{q_i}(u_1a)|} .
$$

It comes that $\psi (\tau^{m_i}(b))$ (resp. $\phi (\sigma^{q_i}(a))$) has an occurrence in $\phi (\sigma^{n_i}(a'))$ (resp. $\psi (\tau^{p_i}(b'))$).

\medskip

To obtain a contradiction it suffices to prove that there exists $j\geq i_0$ such that 
$$
\gap(\sigma^{n_j} (a'))/ |\tau^{m_j} (b)|\leq \frac{1}{2} \ \ {\rm or} \ \ \gap(\tau^{p_j} (b'))/ |\sigma^{q_j} (a)|\leq \frac{1}{2}.
$$

We will consider several cases. Before we define $K$ to be the maximum
of the set 
$$
\left\{  K^{'},
2 \frac{\log \beta}{\log \alpha},
2 \frac{\log \alpha}{\log \beta},
\frac{4\lambda_{\tau} (v_1)}{2 \lambda_{\sigma}(u_1) + \lambda_{\sigma}(a')}
,
\frac{4\lambda_{\sigma} (u_1)}{2 \lambda_{\tau}(v_1) + \lambda_{\tau}(b')}
\right\} .
$$
We remark that $K\geq 2$. 
There exists $j_0$ such that for all $i\geq j_0$ the quantities
$$
\frac{n_i}{m_i} 
\ , \ \ 
\frac{p_i}{q_i} 
\ , \ \ 
\frac{n_i^d \alpha^{n_i}}{m_i^e\beta^{m_i}}
\ , \ \
\frac{p_i^e \beta^{p_i}}{q_i^d\alpha^{q_i}}
\ , \ \
\frac{\mu (a) q_i^{d'}{\alpha^{'}}^{q_i}}{|\sigma^{q_i} (a)|}
\ , \ \
\frac{\mu (b) m_i^{e'}{\beta^{'}}^{m_i}}{|\tau^{m_i} (b)|}
\ \ {\rm and} \ \
\frac{\gap(\sigma^{n_i} (a'))}{ n_i^{d'} {\alpha^{'}}^{n_i}}
$$
are less than $K$. Let $i\geq j_0$. To find $j$ we will consider five cases.

\medskip 

{\bf First case: $\frac{\log (\alpha)}{\log (\beta)} < \frac{\log
(\alpha^{'})}{\log (\beta^{'})}$}. As $\beta^{'}>1$ we have
$$
\gap(\tau^{p_i} (b')/ |\sigma^{q_i} (a)| 
\leq 
\frac{K p_i^{e^{'}} {\beta^{'}}^{p_i}}{\mu (a) q_i^{d'} {\alpha^{'}}^{q_i}} \frac{\mu (a) q_i^{d^{'}} {\alpha^{'}}^{q_i}}{|\sigma^{q_i} (a)|}
$$
$$
\leq
\frac{K^2}{\mu (a)} \frac{p_i^{e'}}{q_i^{d'}} 
\exp
\left\{
\left( 
\frac{p_i}{q_i} - \frac{\log \alpha^{'}}{\log {\beta^{'}}}
\right)
q_i \log \beta^{'}
\right\}
,
$$
which tends to $0$ when $i$ tends to $\infty$ (this comes from (\ref{nu})).

\medskip
 
{\bf Second case: $\frac{\log (\alpha^{'})}{\log (\beta^{'})} < \frac{\log (\alpha)}{\log (\beta)}$}. As in the first case we obtain
$$
\lim_{i\rightarrow +\infty} \gap(\sigma^{n_i} (a'))/ |\tau^{m_i} (b)| = 0.
$$

{\bf Third case: $\frac{\log (\alpha^{'})}{\log (\alpha)} = \frac{\log (\beta^{'})}{\log (\beta)}$ and $(e^{'} - d^{'})\log \beta < (e-d) \log \beta^{'}$}. We have
$$
\gap(\tau^{p_i} (b'))/ |\sigma^{q_i} (a)| 
\leq 
\frac{K^2}{\mu (a)} \frac{p_i^{e'}}{q_i^{d'}} 
\frac{ {\beta^{'}}^{p_i}}{ {\alpha^{'}}^{q_i}}
=
\frac{K^2}{\mu (a)}\frac{p_i^{e'}}{q_i^{d'}} 
\left(
\frac{ {\beta}^{p_i}}{ {\alpha}^{q_i}}
\right) ^{\frac{\log \beta^{'}}{\log \beta}}
$$
$$
= \frac{K^2}{\mu (a)}
\frac{p_i^{e'}}{q_i^{d'}} 
\left(
\frac{q_i^{d}}{p_i^{e}}
\right) ^{\frac{\log \beta^{'}}{\log \beta}} 
\left(
\frac{p_i^{e}\beta^{p_i}}{q_i^{d}\alpha^{q_i}}
\right) ^\frac{\log \beta^{'}}{\log \beta}
$$
$$
\leq
\frac{K^2}{\mu (a)} 
\left(
\frac{p_i}{q_i}
\right) ^{e^{'}-e\frac{\log \beta^{'}}{\log \beta}}
K^{\frac{\log \beta^{'}}{\log \beta}}
q_i^{( e^{'} - d^{'}) - (e - d)\frac{\log \beta^{'}}{\log \beta}}
$$
$$\leq
\frac{K^2}{\mu (a)} K^{e^{'}+(1-e)\frac{\log \beta^{'}}{\log \beta}}
q_i^{( e^{'} - d^{'} ) - (e - d)\frac{\log \beta^{'}}{\log \beta}},
$$
which tends to 0 when $i$ tends to $\infty$.

\medskip

{\bf Fourth case:  $\frac{\log (\alpha^{'})}{\log (\alpha)} = \frac{\log
(\beta^{'})}{\log (\beta)}$ and $(e^{'} - d^{'})\log \beta > (e-d) \log \beta^{'}$}.  As in the previous case we obtain
$$
\lim_{i\rightarrow +\infty} \gap(\sigma^{n_i} (a'))/ |\tau^{m_i} (b)| = 0.
$$

\medskip

{\bf Fifth case:  $\frac{\log (\alpha^{'})}{\log  (\alpha)} = \frac{\log
(\beta^{'})}{\log (\beta)}$ and $(e^{'} - d^{'})\log \beta = (e-d) \log
\beta^{'}$}. From (\ref{correction}), (\ref{multind}) and (\ref{nu}) we obtain
for all large enough $i$
$$
\gap(\tau^{p_i} (b'))/ |\sigma^{q_i} (a)| 
\leq
\frac{K'}{\mu (a)} \frac{p_i^{e^{'}} {\beta^{'}}^{p_i}}{ q_i^{d'} {\alpha^{'}}^{q_i}} \frac{\mu (a) q_i^{d^{'}} {\alpha^{'}}^{q_i}}{|\sigma^{q_i} (a)|}
$$
$$
\leq
\frac{K'}{\mu (a)}
\left( 
\frac{p_i^{e} {\beta}^{p_i}}{ q_i^{d} {\alpha}^{q_i}} 
\right)^{\frac{\log \alpha' }{ \log \alpha}}
\left(
\frac{q_i}{p_i}
\right)^{e\frac{\log \beta'}{\log \beta} - e'}
\frac{\mu (a) q_i^{d^{'}} {\alpha^{'}}^{q_i}}{|\sigma^{q_i} (a)|}
\leq
\frac{1}{2}.
$$
This ends the proof.\hfill $\Box$

\begin{coro}
\label{wordgap}
The words having infinitely many occurrences in $\x$ appear in $\x$ with bounded gaps.
\end{coro}
{\bf Proof.} 
Let $u$ be a word having infinitely many occurrences in $\x$. We set $|u| = n$.
 To prove that $u$ appears with bounded gaps in $\x$ it suffices to prove that the letter 1 appears with bounded gaps in the sequence $\t \in \{ 0,1 \}^{\NN}$ defined by
$$
\t_i = 1 \ \ {\rm if} \ \ \x_{[i,i+n-1]} = u
$$
and $0$ otherwise. In the sequel we prove that $\t$ is $\alpha$ and $\beta$-substitutive.

The sequence $\y^{(n)} = (( \y_i \cdots \y_{i+n-1}) ; i\in \NN )$ is a fixed  point of the substitution $\sigma_n : A_n \rightarrow A_n^*$ where $A_n$ is the alphabet $A^n $, defined for all $(a_1\cdots a_n)$ in $A_n$ by
$$
\sigma_n ((a_1\cdots a_n)) = (b_1\cdots b_n)(b_2\cdots b_{n+1})\cdots (b_{|\sigma (a_1)|}\cdots b_{|\sigma (a_1)|+n-1})
$$
where $\sigma(a_1\cdots a_n) = b_1\cdots b_k$ (for more details see
Section V.4 in \cite{Qu} for example).

Let $\rho : A_n \rightarrow A^{*}$ be the letter to letter morphism defined by  $\rho ( ( b_1\cdots b_n )) = b_1$ for all $(b_1\cdots b_n)\in A_n$. We have $\rho \circ \sigma_n = \sigma \circ \rho$, and then $M_{\rho}M_{\sigma_n} = M_{\sigma} M_{\rho}$. Consequently the dominant eigenvalue of $\sigma_n$ is $\alpha$ and $\y^{(n)}$ is $\alpha$-substitutive.

Let $f : A_n \rightarrow \{ 0,1 \}$ be the letter to letter morphism
 defined by

\medskip

\centerline
{
$f ( ( b_1\cdots b_n )) = 1$ if $b_1\cdots b_n = u$ and $0$ otherwise.
}

\medskip

It is easy to see that $f( \y^{(n)} ) = \t$ hence $\t$ is $\alpha$-substitutive.

In the same way we show that $\t$ is $\beta$-substitutive and Theorem \ref{bgaps} concludes the proof.
\hfill $\Box$

\section{Proof of Theorem 1.}

\label{proof}

\subsection{Decomposition of a substitution into sub-substi\-tu\-tions.}

The following proposition is a consequence of Paragraph 4.4 and Proposition 4.5.6 in \cite{LM}.
\begin{prop}
\label{decomprim}
Let $M=(m_{i,j})_{i,j\in A}$ be a matrix with non-negative coefficients
 and no column equal to 0. There exist three positive integers $p\not = 0$, $q$, $l$, where $q\leq l-1$, and a partition $\{ A_i ; 1\leq i\leq l \}$ of $A$ such that the matrix $M^p$ is equal to 
\begin{equation}
\label{formefin}
\bordermatrix{        & A_1     & A_2         & \cdots & A_q       & A_{q+1} & A_{q+2} & \cdots & A_{l}  \cr 
                  A_1     & M_1     & 0           & \cdots & 0         & 0       & 0       & \cdots & 0      \cr
                  A_2     & M_{1,2} & M_2         & \cdots & 0         & 0       & 0       & \cdots & 0      \cr
                  \vdots  & \vdots  & \vdots      & \ddots & \vdots    & \vdots  & \vdots  & \vdots & \vdots \cr
                  A_q     & M_{1,q} & M_{2,q}     & \cdots & M_q       & 0       & 0       & \cdots & 0      \cr
                  A_{q+1} & M_{1,q+1} & M_{2,q+1} & \cdots & M_{q,q+1} & M_{q+1} & 0       & \cdots & 0      \cr
                  A_{q+2} & M_{1,q+2} & M_{2,q+2} & \cdots & M_{q,q+2} & 0       & M_{q+2} & \cdots & 0      \cr
                  \vdots  & \vdots    & \vdots    & \ddots & \vdots    & \vdots  & \vdots  & \ddots & \vdots \cr
                  A_l     & M_{1,l}   & M_{2,l}   & \cdots & M_{q,l}   & 0       & 0       & \cdots & M_l    \cr}, 
\end{equation}
where the matrices $M_i$, $1\leq i\leq q$ (resp. $q+1\leq i\leq l$) , are primitive or equal to zero (resp. primitive), and such that for all $1\leq i\leq q$ there exists $i+1 \leq j \leq l$ such that the matrix $M_{i,j}$ is different from 0.
\end{prop}

In what follows we keep the notations of Proposition \ref{decomprim}.
We will say that $\{ A_i; 1\leq i\leq l \}$ is a {\it primitive component partition of $A$ (with respect to $M$)}. If $i$ belongs to $\{ q+1 , \cdots ,l \}$ we will say that $A_i$ is a {\it principal primitive component of $A$ (with respect to $M$)}.

Let $\tau : A\rightarrow A^{*}$ be a substitution and $M =
(m_{i,j})_{i,j\in A}$ its matrix. Let $i\in \{ q+1,\cdots ,l \}$. We
denote $\tau_i$ the restriction $(\tau^p)_{|A_i} : A_i\rightarrow A^*$
of $\tau^p $ to $A_i$. Because $\tau_i (A_i)$ is included in $A_i^*$ we
can consider that $\tau_i$ is a morphism from $A_i$ to $A_i^*$ whose
matrix is $M_i$. Let $i\in \{ 1,\cdots ,q \}$ such that $M_i$ is  not
equal to 0. Let $\varphi_i$ be the morphism from $A$ to $A_i^*$ defined
by $\varphi (b) = b$ if $b$ belongs to $A_i$ and the empty word
otherwise. Let us consider the map $\tau_i : A_i \rightarrow A^*$ defined by $\tau_i (b) = \varphi_i (\tau^p (b))$ for all $a\in A_i$. We remark as previously that $\tau_i (A_i)$ is included in $ A_i^*$, consequently $\tau_i$ defines a morphism from $A_i$ to $A_i^*$ whose matrix is $M_i$.

We will say that the substitution $\tau : A \rightarrow A^*$ satisfies Condition {\bf (C)} if: 
\begin{enumerate}
\item[C1.]
The matrix $M$, itself, is of the type (\ref{formefin}) (i.e., $p=1$);
\item[C2.]
The matrices $M_i$ are equal to 0 or with positive coefficients if $1\leq i\leq q$ and with positive coefficients otherwise;
\item[C3.]
For all matrices $M_i$ different from 0, with $i\in \{ 1,\cdots ,l \}$,
	  there exists $a_i\in A_i$ such that $\tau_i (a_i) = a_iu_i$
	  where $u_i$ is a non-empty word of $A^*$ if $M_i$ is different
	  from the $1\times 1$ matrix $[1]$ and empty otherwise.
\end{enumerate}
From Proposition \ref{decomprim} every substitution $\tau : A \rightarrow A^{*}$ has a power $\tau^k $ satisfying condition {\bf (C)}. The definition of substitutions implies that for all $q+1\leq i\leq l$ we have $M_i \not = [1]$.

Let $\tau : A \rightarrow A^{*}$ be a substitution satisfying condition
{\bf (C)} (we keep the previous notations). For all $1 \leq i \leq l$
such that $M_i$ is different from 0 and the $1\times 1$ matrix $[1]$, the map $\tau_i : A_i \rightarrow A_i^*$ defines a substitution we will call {\it main sub-substitution of $\tau$} if $i\in \{ q+1,\cdots ,l \}$ and {\it non-main sub-substitution of $\tau$} otherwise. Moreover the matrix $M_i$ has positive coefficients which implies that the substitution $\tau_i $ is primitive. We remark that there exists at least one main sub-substitution.

In \cite{Du3} the following results were obtained and will be used in the sequel.

\begin{lemme}
\label{lemmeun}
Let $\x$ be a proper fixed point of the substitution $\sigma$. Let $\sigmabar : \Abar \rightarrow \Abar^*$ be a main sub-substitution
 of $\sigma$. Then for all $n\in \NN$ and all $a\in A$
 the word $\sigma^n (a)$ appears infinitely many times in $\x$.
\end{lemme}
{\bf Proof.}
The proof is left to the reader.
\hfill $\Box $

\medskip

In \cite{Du3} the following result is obtained and will be used in the sequel.

\begin{theo}
\label{factcob} 
Let $\x$ and $\y$ be respectively a primitive $\alpha$-substitutive
 sequence and a primitive $\beta$-substitu\-tive sequence such that $L (\x) = L (\y)$. Suppose that $\alpha$ and $\beta$ are multiplicatively independent, then $\x$ and $\y$ are periodic.
\end{theo}

\subsection{The conjecture for ``good'' substitutions.}

We do not succeed yet to prove the conjecture given in the introduction but we are able to prove it for a very large family of substitutions. Until we prove the whole conjecture we call them ``good'' substitutions. More precisely, let $\sigma : A\rightarrow A^*$ be a substitution whose dominant eigenvalue is $\alpha$. The substitution $\sigma$ is said to be a {\it ``good'' substitution} if there exists a main sub-substitution whose dominant eigenvalue is $\alpha$.

For example primitive substitutions and substitutions of constant length are ``good'' substitutions. Now consider the following substitution
\[
\begin{array}{llll}
\sigma : & \{ a , 0 , 1\} & \rightarrow & \{ a , 0 , 1\} ^* \\
         & a              & \mapsto     & aa0 \\
         & 0              & \mapsto     & 01 \\
         & 1              & \mapsto     & 0 .
\end{array}
\]
Its dominant eigenvalue is $2$ and it has only one main sub-substitution ($0\mapsto 01$, $1\mapsto 0$) which dominant eigenvalue is $(1+\sqrt{5})/2$, hence it is not a ``good'' substitution. 

\begin{theo}
\label{good}
Suppose that we only consider ``good'' substitutions. Then the conjecture is true.
\end{theo}
{\bf Proof.} We take the notations of the first lines of Section \ref{letterbg}.

Let $\sigmabar : \Abar \rightarrow \Abar^*$ be a main sub-substitution
of $\sigma$. The words of $\x$ appearing infinitely many times in $\x$
appear with bounded gaps (Corollary \ref{wordgap}). Hence using Lemma
\ref{lemmeun} we deduce that for all main sub-substitution $\sigmabar $ of $\sigma$ and
 $\taubar$ of $\tau$ we have $\phi (L(\sigmabar )) = \psi ( L( \taubar ))=L$. From Theorem \ref{factcob} it comes that $L$ is periodic, i.e., there exists a word $u$ such that $L = L(u^{\omega})$ where $|u|$ is the least period. There exists an integer $N$ such that all the words of length $|u|$ appear infinitely many times in $\x_N \x_{N+1} \cdots $. We set $\t = \x_N \x_{N+1} \cdots $ and we will prove that $\t$ is periodic and consequently $\x$ will be ultimately periodic.

The word $u$ appears infinitely many times, consequently it appears with bounded gaps. Let $\R_u$ be the set of return words to $u$ (a word $w$ is a return word to $u$ if $wu\in L(\x)$, $u$ is a prefix of $wu$ and $u$ has exactly two occurrences in $wu$). It is finite. There exists an integer $N$ such that all the words $w\in \R_u \cap L(\x_N \x_{N+1} \cdots )$ appear infinitely many times in $\x$. Hence these words appear with bounded gaps in $\x$. We set $\t = \x_N \x_{N+1} \cdots $ and we will prove that $\t$ is periodic and consequently $\x$ will be ultimately periodic. We can suppose that $u$ is a prefix of $\t$. Then $\t$ is a concatenation of return words to $u$. Let $w$ be a return word to $u$. It appears with bounded gaps hence it appears in some $\phi (\sigmabar^n (a) )$ and there exist two words, $p$ and $q$, and an integer $i$ such that $wu = pu^i q$. As $|u|$ is the least period of $L$ it comes that $wu = u^i$. It follows that $\t = u^{\omega}$.

\hfill $\Box$

\bigskip

{\bf The case of fixed points.}

This part is devoted to the proof of Theorem \ref{maintheo} restricted to fixed points. More precisely we prove:
\begin{coro}
\label{fixedp}
Let $\x$ be a fixed point of the substitution $\sigma : A \rightarrow A^*$ whose dominant eigenvalue is $\alpha$. Suppose that $\x$ is also a fixed point of the substitution $\tau : A \rightarrow A^*$ whose dominant eigenvalue is $\beta$.  Suppose that $\alpha$ and $\beta$ are multiplicatively independent. Then  $\x$ is ultimately periodic.
\end{coro}

{\bf Proof.} The letters appearing infinitely often in $\x$ appear with
bounded gaps (Proposition \ref{bgaps}). Let $\sigmabar : \Abar
\rightarrow \Abar$ be a main sub-substitution of $\sigma$. Let $a\in \Abar$. Suppose that there exists a
letter $b$, appearing infinitely many times in $\x$, which does not
belong to $\Abar$. Then the word $\sigma^n (a)$ does not contain $b$ and
$b$ could not appear with bounded gaps. Consequently there exists only
one main sub-substitution and the letters which appear with bounded gaps
belong to $\Abar$. It comes that $\sigma$ is a ``good'' substitution. In
the same way $\tau$ is a good substitution. Theorem \ref{good} concludes
the proof. 
\hfill $\Box$

\end{document}